\IfFileExists{pdfmanagement.sty}{%
  \RequirePackage{pdfmanagement}
  \SetKeys[document/metadata]{uncompress}
}{%
  \ifdefined\pdfcompresslevel
    \pdfcompresslevel=0\relax
  \fi
  \ifdefined\pdfobjcompresslevel
    \pdfobjcompresslevel=0\relax
  \fi
}
\documentclass[11pt]{amsart}

\usepackage{etoolbox}

\makeatletter
\patchcmd{\@setauthors}
  {\MakeUppercase{\authors}}
  {\authors}
  {}
  {\PackageWarning{authorcase}{Author patch failed}}
\makeatother

 \usepackage{amsfonts,graphics,amsmath,amsthm,amsfonts,amscd,amssymb,amsmath,latexsym,multicol,
 mathrsfs}
\usepackage{epsfig,url}
\usepackage{flafter}
\usepackage{fancyhdr}
\usepackage{hyperref}
\usepackage[dvipsnames]{xcolor}
\usepackage{hyperref}
\hypersetup{
	colorlinks=true,
	linkcolor=purple,
	citecolor=ForestGreen,
	filecolor=cyan,
	urlcolor=purple
}

\usepackage{graphicx}
\usepackage{xcolor}
\usepackage{float}
\usepackage{bm}
\usepackage{enumitem}
\usepackage{multirow}
\usepackage{scalerel,stackengine}
\usepackage{stmaryrd}
\usepackage{datetime}
\usepackage[title]{appendix}
\usepackage{url}
\usepackage{tabulary}
\usepackage{booktabs}
\usepackage{tikz}
\usetikzlibrary{positioning}
\usepackage{verbatim}

\usepackage[utf8]{inputenc}
\usepackage[english]{babel}
\usepackage{csquotes}
\usepackage{comment}
\usepackage[
backend=biber,
style=alphabetic,
sorting=anyt,
doi=false,isbn=false,url=false,
maxalphanames=99,maxbibnames=99,backref=true,
]{biblatex}
\DefineBibliographyStrings{english}{
  backrefpage  = {\hspace{-0.35em}},
  backrefpages = {\hspace{-0.35em}},
}

\DeclareFieldFormat{bracketswithperiod}{\mkbibbrackets{#1}}

\renewbibmacro*{pageref}{%
  \iflistundef{pageref}
    {}
    {\printtext[bracketswithperiod]{%
       \ifnumgreater{\value{pageref}}{1}
         {\bibstring{backrefpages}\ppspace}
         {\bibstring{backrefpage}\ppspace}%
       \printlist[pageref][-\value{listtotal}]{pageref}}%
     }}

 \theoremstyle{plain}
 \newtheorem{theorem}[subsection]{Theorem}
 \newtheorem{lemma}[subsection]{Lemma}

 \theoremstyle{definition}

 \theoremstyle{remark}

 \numberwithin{equation}{subsection}

 \theoremstyle{plain} 
\newcommand{\thistheoremname}{}
\newtheorem{genericthm}[subsection]{\thistheoremname}

\newcommand{\abs}[1]{\left\vert#1\right\vert}
\newcommand{\set}[1]{\left\{#1\right\}}

\newcommand{\spec}{\operatorname{Spec}}

\newcommand{\codim}{\operatorname{codim}}

\newcommand{\chara}{\operatorname{char}}
\newcommand{\projection}{\operatorname{pr}}

\newcommand{\Supp}{\operatorname{Supp}}

\newcommand{\sdeg}{\operatorname{sdeg}}

\newcommand{\centre}{\operatorname{centre}}

\newcommand{\Ex}{\operatorname{Ex}}

\newcommand{\coeff}{\operatorname{coeff}}
\newcommand{\nor}{\operatorname{nor}}

\newcommand{\bb}{\operatorname{b}}

\newcommand{\N}{\operatorname{N}}
\newcommand{\Gal}{\operatorname{Gal}}

\newcommand{\NE}{\operatorname{\overline{NE}}}
\newcommand{\Nef}{\operatorname{Nef}}

\newcommand{\univ}{\operatorname{univ}}

\newcommand{\sym}{\operatorname{Sym}}

\newcommand{\mni}{\medskip\noindent}

\newcommand{\mbb}{\mathbb}
\newcommand{\QQ}{\mbb{Q}}
\newcommand{\NN}{\mbb{N}}
\newcommand{\ZZ}{\mbb{Z}}
\newcommand{\CC}{\mbb{C}}
\newcommand{\RR}{\mbb{R}}

\newcommand{\PP}{\mbb{P}}

\newcommand{\mc}{\mathcal}
\newcommand{\mcA}{\mc{A}}
\newcommand{\mcB}{\mc{B}}

\newcommand{\mcO}{\mc{O}}

\newcommand{\mf}{\mathfrak}

\newcommand{\mb}{\mathbf}

\newcommand{\wt}{\widetilde}
\newcommand{\ol}{\overline}
\newcommand{\ul}{\underline}
\newcommand{\bir}{\dashrightarrow}

\usepackage{graphicx}
\usepackage{xcolor}
\usepackage{float}
\usepackage{amssymb}
\usepackage{amsmath}
\usepackage{mathrsfs}
\usepackage{bm}
\usepackage[all,cmtip]{xy}
\usepackage{enumitem}
\usepackage[utf8]{inputenc}
\usepackage[english]{babel}
\usepackage{multirow}
\usepackage{mathtools}
\usepackage{scalerel,stackengine}
\usepackage{stmaryrd}
\usepackage{datetime}
\usepackage[title]{appendix}
\usepackage{url}
\usepackage{tabulary}
\usepackage{booktabs}
\usepackage{tikz}
\usetikzlibrary{positioning}
\usepackage{verbatim}
\usepackage{braket}

\DeclarePairedDelimiter\rddown{\lfloor}{\rfloor}
\DeclareFieldFormat{postnote}{#1}
\DeclareFieldFormat{multipostnote}{#1}

\makeatletter
\def\moverlay{\mathpalette\mov@rlay}
\def\mov@rlay#1#2{\leavevmode\vtop{%
   \baselineskip\z@skip \lineskiplimit-\maxdimen
   \ialign{\hfil$\m@th#1##$\hfil\cr#2\crcr}}}
\newcommand{\charfusion}[3][\mathord]{
    #1{\ifx#1\mathop\vphantom{#2}\fi
        \mathpalette\mov@rlay{#2\cr#3}
      }
    \ifx#1\mathop\expandafter\displaylimits\fi}
\makeatother

\makeatletter
\newcommand*{\da@rightarrow}{\mathchar"0\hexnumber@\symAMSa 4B }
\newcommand*{\da@leftarrow}{\mathchar"0\hexnumber@\symAMSa 4C }
\newcommand*{\xdashrightarrow}[2][]{%
  \mathrel{%
    \mathpalette{\da@xarrow{#1}{#2}{}\da@rightarrow{\,}{}}{}%
  }%
}
\newcommand{\xdashleftarrow}[2][]{%
  \mathrel{%
    \mathpalette{\da@xarrow{#1}{#2}\da@leftarrow{}{}{\,}}{}%
  }%
}
\newcommand*{\da@xarrow}[7]{%
  \sbox0{$\ifx#7\scriptstyle\scriptscriptstyle\else\scriptstyle\fi#5#1#6\m@th$}%
  \sbox2{$\ifx#7\scriptstyle\scriptscriptstyle\else\scriptstyle\fi#5#2#6\m@th$}%
  \sbox4{$#7\dabar@\m@th$}%
  \dimen@=\wd0 %
  \ifdim\wd2 >\dimen@
    \dimen@=\wd2 %
  \fi
  \count@=2 %
  \def\da@bars{\dabar@\dabar@}%
  \@whiledim\count@\wd4<\dimen@\do{%
    \advance\count@\@ne
    \expandafter\def\expandafter\da@bars\expandafter{%
      \da@bars
      \dabar@ 
    }%
  }%
  \mathrel{#3}%
  \mathrel{%
    \mathop{\da@bars}\limits
    \ifx\\#1\\%
    \else
      _{\copy0}%
    \fi
    \ifx\\#2\\%
    \else
      ^{\copy2}%
    \fi
  }%
  \mathrel{#4}%
}
\makeatother

\title{\large S\MakeLowercase{tein degree on birational log} C\MakeLowercase{alabi-}Y\MakeLowercase{au fibrations}}
\author[Xintong Jiang and Santai Qu]{%
  \large Xintong Jiang and Santai Qu
}
\date{\today}

\begin{document}

\begin{abstract}
    We prove uniform bounds for the Stein degrees of boundary divisors
    on generalised log Calabi-Yau fibrations of Fano type over arbitrary
    fields of characteristic zero.  The main ingredient is
    the birational case for which we provide two approaches.
\end{abstract}

\maketitle

\tableofcontents

\addtocontents{toc}{\protect\setcounter{tocdepth}{1}}


\section{\bf Introduction}

\mni
{\textbf{\sffamily{Stein degree on log Calabi-Yau fibrations.}}}
The fibres of a birational contraction between normal schemes are geometrically connected, but
an individual exceptional divisor can carry a nontrivial finite cover of
its centre.  In the presence
of a log Calabi-Yau structure, the degree of this cover connects the singularities
of the ambient pair with the geometry and arithmetic of its boundary
divisors.  The purpose of this paper is to establish uniform bounds for this degree.

Let $S \to T$ be a projective surjective morphism between integral schemes,
and let
\[ S \to V \to T \]
be its Stein factorisation.  The \emph{Stein degree} is defined to be
\[ \sdeg (S/T) \coloneqq \deg(V/T). \]
For a prime divisor $S$ on an integral scheme $X$ admitting a contraction $X \to Z$,
we consider the Stein degree of its normalisation $S^{\nor}$
over its image $T \subseteq Z$.  In characteristic zero,
$\sdeg (S^{\nor}/T)$ counts the irreducible components of the geometric generic fibre of $S \to T$. 
Equivalently, it is the degree of the relative algebraic closure of $K(T)$ in $K(S)$.
This interpretation makes Stein degree a natural invariant for studying
how a divisor splits over its image.

The boundedness of Stein degree arose in Birkar's work on moduli of
algebraic varieties \cite{B-moduli}. For a log Calabi-Yau fibration $(X,B)\to Z$
over an algebraically closed field of characteristic zero, \cite{B-moduli}
proved that the Stein degrees $\sdeg (I/Z)$ of horizontal
log canonical centres $I$ are bounded in terms of $\dim X$.  This result
is a crucial ingredient in the boundedness of semi-log canonical stable
minimal models.  Birkar further conjectured that a similar bound should
hold for any horizontal boundary component whose coefficient is
bounded below by a fixed positive number, thereby extending the
question beyond log canonical centres; see
\cite[Theorem~1.15, Conjecture~1.16]{B-moduli}.

Birkar and the second author of this paper 
subsequently proved \cite[Conjecture~1.16]{B-moduli} in the more general
setting of generalised pairs \cite{birkar2025steindegreelogcalabiyau}.  
Over an algebraically closed field of
characteristic zero, \cite{birkar2025steindegreelogcalabiyau} showed that for a generalised log
Calabi-Yau fibration of relative dimension $d$, the Stein degree
\[ \sdeg (S^{\nor}/Z) \]
is bounded from above in terms of $d$ and a positive lower bound $t$ for
the coefficient of a horizontal boundary irreducible component $S$. 
\cite{birkar2025steindegreelogcalabiyau} also
established a vertical counterpart: if the contraction is of Fano type
and the image $T$ of $S$ is a divisor on $Z$, then
$\sdeg (S^{\nor}/T)$ admits an upper bound depending only
on $d,t$; see \cite[Theorems~1.3 and~1.4]{birkar2025steindegreelogcalabiyau}.
Note that no Fano type assumption is required for the horizontal result. 

\mni
{\textbf{\sffamily{Stein degree of exceptional divisors.}}}
The two boundedness statements of \cite{birkar2025steindegreelogcalabiyau} were extended to arbitrary fields
of characteristic zero in \cite[Theorems~1.1 and~1.2]{birkar2026steindegreefieldscharacteristic}, with
constants independent of the ground field.  In the absolute setting,
a uniform arithmetic bound for boundary irreducible components of ordinary
projective log Calabi-Yau pairs with rational boundaries was also
independently announced in \cite[Theorem~1.1]{OAI:Arithmetic-Stein-degree-bounds-for-log-Calabi-Yau-pairs-September-25-2026}.

The remaining case concerns boundary divisors whose images have
codimension at least two.  Such divisors do not meet the general fibre
of the original contraction, and their behaviour is not covered by
the results in \cite{birkar2025steindegreelogcalabiyau, birkar2026steindegreefieldscharacteristic}. 
The issue is
already present for birational contractions: although the morphism
itself has degree one, the normalisation of an exceptional divisor can
have large Stein degree over its centre.  We prove that this degree
is nevertheless uniformly bounded 
if the birational contraction is of Fano type.

Throughout the paper, $k$ is an arbitrary field of characteristic zero.
Recall that a \emph{generalised log Calabi-Yau fibration} $(X, B + \mb{M}) \to Z$ over $k$ 
consists of a generalised log canonical 
generalised pair $(X, B + \mb{M})$ over $k$ and a contraction between $k$-varieties $X\to Z$ 
such that $K_X + B + \mb{M}_X \sim_{\RR}0/Z$.
Moreover, a contraction $X \to Z$ is of \emph{Fano type} if there exists an
effective $\QQ$-divisor $\Gamma$ such that $(X,\Gamma)$ is klt
and $-(K_X+\Gamma)$ is ample$/Z$.

Our first main result treats birational contractions.

\begin{theorem}[$\mc{B}_{\univ}(d,t)$]\label{dim-d-t}
    Let $d\in \NN$ and let $t\in (0,1]$.
    There exists an $N_b(d,t) \in \NN$ depending only on $d,t$ such that
    the following holds.

    Let $k$ be a field of characteristic zero.
	Let $(X, B + \mb{M})\to Z$ be a generalised log Calabi-Yau fibration over $k$,
	where $X\to Z$ is a birational contraction of Fano type between $k$-varieties, and $\dim X = d$.
    Let $S$ be an irreducible component of $B$ 
	whose coefficient in $B$ is $\ge t$.  Then 
    \[ \sdeg (S^{\nor}/T) \le N_b(d,t), \]
    where $T$ is the image of $S$ in $Z$, and $S^{\nor}$ is the normalisation of $S$.
\end{theorem}

When $S$ is not exceptional over $Z$, the morphism $S \to T$ is
birational and the Stein degree is one.  Thus, the substance of
Theorem~\ref{dim-d-t} lies in its control of
exceptional divisors. 


\mni
{\textbf{\sffamily{Stein degree of vertical divisors.}}}
The birational statement Theorem~\ref{dim-d-t} also supplies the
missing ingredient for treating boundary irreducible components on general
Fano type fibrations, without any restriction on the codimension of their images.

\begin{theorem}[$\mbb{V}_{\univ}(d,t)$]\label{thm-vertical-introduction}
    Let $d\in \NN$ and let $t\in (0,1]$.  There exists an $N_{\mbb V}(d,t)\in \NN$     
    depending only on $d,t$ such that the following holds.

    Let $k$ be a field of characteristic zero.  Let $(X, B + \mb M)\to Z$
    be a Fano type generalised log Calabi-Yau fibration over $k$ of dimension $d$, that is,
    \begin{itemize}
        \item $X, Z$ are normal $k$-varieties with $\dim X = d$,
        \item $(X, B + \mb M)\to Z$ is a generalised log Calabi-Yau fibration over $k$, and
        \item $X\to Z$ is a contraction of Fano type between $k$-varieties.
    \end{itemize}
    Let $S$ be an irreducible component of $B$ whose coefficient in $B$ is $\ge t$.  Then
    \[ \sdeg (S^{\nor}/T) \le N_{\mbb V}(d,t), \]
    where $T$ is the image of $S$ in $Z$, and $S^{\nor}$ is the normalisation of $S$.
\end{theorem}

Theorem~\ref{thm-vertical-introduction} applies to both horizontal
and vertical boundary irreducible components.  Its new content is the uniform
bound for vertical irreducible components with higher-codimensional images.
Note that the constants remain independent of the ground field.


\mni
{\textbf{\sffamily{Outline of the proofs.}}}
We give two proofs of Theorem~\ref{dim-d-t}. In both arguments,
the Fano type assumption allows us to run suitable minimal model
programmes and reduce to ordinary log Calabi-Yau pairs with
rational boundaries.  Without this assumption, Stein degrees of
vertical divisors need not be bounded; see
\cite{birkar-non-fano-stein}.

The first proof uses a result of 
\cite{Hogadi-Xu-09}, which provides geometrically irreducible
subschemes in fibres of birational contractions onto varieties
of klt type.  The second proof proceeds by simultaneous induction
on dimension for Theorems~\ref{dim-d-t}
and~\ref{thm-vertical-introduction}.  The induction combines a
reduction of the birational statement to the general statement
in lower dimensions with a deduction of
Theorem~\ref{thm-vertical-introduction} from
Theorem~\ref{dim-d-t} and
\cite[Theorem~1.2]{birkar2026steindegreefieldscharacteristic}.

Both approaches rely on
\cite[Theorems~1.1 and~1.2]{birkar2026steindegreefieldscharacteristic}.
Thus, the present paper builds on these results and does not
provide alternative proofs of them.


\mni
{\textbf{\sffamily{Comparison with \cite{birkar2026steindegreefieldscharacteristic}.}}}
Keep the notation of Theorem~\ref{thm-vertical-introduction}. 
When $T=Z$, the required bound follows from
\cite[Theorem~1.1]{birkar2026steindegreefieldscharacteristic},
without the Fano type assumption.  When $T$ is a divisor on~$Z$,
it follows from \cite[Theorem~1.2]{birkar2026steindegreefieldscharacteristic}.
In these two cases, the bounds depend only on the relative
dimension of $X\to Z$ and the coefficient lower bound~$t$.
Theorem~\ref{thm-vertical-introduction} allows $T$ to have arbitrary codimension, 
while the bound of Stein degree depends on the total dimension of~$X$.

The examples in \S\ref{s-examples} show that this distinction is necessary.
In \S\ref{exa-dim}, we construct birational log Calabi-Yau fibrations of
Fano type with a boundary irreducible component of coefficient one whose
Stein degree over its image is arbitrarily large as the
dimension increases.  Thus, the dependence on dimension in
Theorem~\ref{dim-d-t} cannot be omitted.  Taking products with
$\PP^1$ gives the examples in \S\ref{exa-rel-dim},
which have relative dimension one and the same Stein degrees.
Consequently, for vertical boundary irreducible components whose images
have codimension at least two, the relative dimension and
the coefficient lower bound alone do not suffice to bound
the Stein degree.  These examples already exist over $\CC$.


\mni
{\textbf{\sffamily{Organisation of the paper.}}}
In \S\ref{s-preliminaries}, we recall the necessary preliminaries and establish
auxiliary results on Fano type contractions, \'etale
neighbourhoods, and plt blow-ups.  In \S\ref{s-first-pf}, we give the first
proof of Theorem~\ref{dim-d-t}.  In \S\ref{s-second-pf}, we develop the inductive proof,
establish the reduction from Theorem~\ref{thm-vertical-introduction} to the birational
case, and give the second proof of Theorem~\ref{dim-d-t}. 
Finally, \S\ref{s-examples} contains examples showing the necessity
of dimension dependence and the failure of bounds depending
only on relative dimension and the coefficient lower bound.


\mni
{\textbf{\sffamily{Acknowledgements.}}}
We thank Professor Caucher Birkar for encouraging us to pursue this project
and for many helpful discussions.
This work was completed during ICBS in August 2026.
The manuscript was prepared after
\cite{birkar2026steindegreefieldscharacteristic} became available
so that we could include a reference to it.


\mni
{\textbf{\sffamily{AI use disclosure.}}}
The main ideas and overall direction of this paper were developed by the authors.
During the preparation of the manuscript, AI tools (ChatGPT models) were used
as technical assistants to help check arguments, identify relevant references,
and improve the exposition.
The examples in \S\ref{s-examples} were suggested by ChatGPT-6 Astra
and worked out in detail by the authors.
All mathematical statements, proofs, and conclusions were verified by the authors.


\section{\bf Preliminaries}\label{s-preliminaries}

Throughout, we work over fields of characteristic zero.
We adopt the notation and conventions of \cite{birkar2026steindegreefieldscharacteristic}.
We indicate below the background used in this paper and
establish the additional auxiliary lemmas needed for our arguments.


\subsection{Schemes, contractions, and Stein degree}

For schemes and varieties over arbitrary fields, contractions,
generic fibres, and Stein factorisation, we refer to
\cite[\S\S2.1-2.4]{birkar2026steindegreefieldscharacteristic}.
In particular, a \emph{variety} is a quasi-projective, integral $k$-scheme
that is geometrically integral;
when emphasising the quasi-projectivity, we also call it a \emph{quasi-projective variety}.
We use the interpretation of Stein degree in terms of
geometric fibres, its birational invariance after
normalisation, and the comparison results in
\cite[Lemmas~2.5, 2.7, 2.9, and~7.1]{birkar2026steindegreefieldscharacteristic}.


\subsection{Pairs and singularities}

For ordinary pairs and their singularities, we follow
\cite[\S2.27]{birkar2026steindegreefieldscharacteristic}.
For $\bb$-divisors, generalised pairs, and
generalised dlt models, we use \cite[\S\S2.30-2.31 and~3.10]{birkar2026steindegreefieldscharacteristic}.
The conventions for flat pullbacks and the compatibility
of singularities with field extensions and \'etale base
change are those of \cite[\S\S2.13 and~2.34-2.35]{birkar2026steindegreefieldscharacteristic}.


\subsection{Fano type and the MMP}

For Fano type contractions and their behaviour under
field extensions and birational operations, we refer to
\cite[\S2.46]{birkar2026steindegreefieldscharacteristic}.
The minimal model programmes (MMPs) and extraction of divisors used in this paper
are supplied by \cite[Theorem~3.8 and Lemmas~3.9 and~3.14]{birkar2026steindegreefieldscharacteristic}.

\begin{lemma}\label{birational-fano-type}
    Let $k$ be a field of characteristic zero.  Let $f\colon X\to Y$ be a contraction$/Z$ of
    $k$-varieties $X,Y$ that are projective over a $k$-variety $Z$.  Assume that
    $X$ is of Fano type over $Z$.  Then $Y$ is of Fano type over $Z$.
\end{lemma}

\begin{proof}
    By \cite[Lemma~2.47]{birkar2026steindegreefieldscharacteristic}, 
    we can assume that $k$ is algebraically closed. 
    Pick an effective $\QQ$-divisor $\Gamma$ such that $(X,\Gamma)$
    is klt, $\Gamma$ is big$/Z$, and $K_X+\Gamma\sim_{\QQ}0/Z$.
    By \cite[Lemma~3.14]{birkar2026steindegreefieldscharacteristic}, 
    we can take a small $\QQ$-factorialisation of $(X,\Gamma)$. 
    By \cite[Lemma~2.50]{birkar2026steindegreefieldscharacteristic},
    this preserves Fano type over $Z$.  Thus, we can assume that $X$ is $\QQ$-factorial.
    
    If $\dim(Y/Z)>0$, the conclusion follows from \cite[Lemma~2.52]{birkar2026steindegreefieldscharacteristic}. 
    We may therefore assume that $Y\to Z$ is birational. 
    In this case, $\Gamma$ is big$/Y$ and $K_X+\Gamma\sim_{\QQ}0/Y$.
    Pick rational numbers $0<\epsilon'<\epsilon$ such that $(X,\Gamma)$ is $\epsilon$-lc. 
    Then \cite[Lemma~3.3]{birkar-sing-base-fano} gives projective compactifications $X\subset X'$ and $Y\subset Y'$ 
    and a contraction $f'\colon X'\to Y'$ extending $f$, with a rational boundary $\Gamma'$ such that 
    \begin{itemize}
        \item $(X',\Gamma')$ is $\epsilon'$-lc, 
        \item $(K_{X'}+\Gamma')|_X = K_X+\Gamma$, and $K_{X'}+\Gamma'\sim_{\QQ}0/Y'$.
    \end{itemize}
    By \cite[Theorem~0.2]{Ambro-lc-fibration}, there exists an effective $\QQ$-divisor $\Gamma_{Y'}$ 
    such that $(Y',\Gamma_{Y'})$ is klt and
    \[ K_{X'}+\Gamma'\sim_{\QQ} (f')^*(K_{Y'}+\Gamma_{Y'}). \]
    Let $\Gamma_Y \coloneqq \Gamma_{Y'}|_Y$.  Then $(Y,\Gamma_Y)$ is klt and $K_Y+\Gamma_Y\sim_{\QQ}0/Z$.
    Since $Y\to Z$ is birational, $\Gamma_Y$ is big$/Z$.  Thus, $Y\to Z$ is of Fano type.
\end{proof}


\begin{lemma}\label{fano-type-polyhedral}
    Let $k$ be a field of characteristic zero.
    Let $X\to Z$ be a contraction of Fano type between $k$-varieties.
    Then $\NE(X/Z)$ is a rational polyhedral cone.
\end{lemma}

\begin{proof}
    As $X\to Z$ is of Fano type, there exists a $\QQ$-divisor $\Gamma$ such that $(X, \Gamma)$
    is klt and $A\coloneqq -(K_X+\Gamma)$ is ample$/Z$.  
    Denote by $\mcA(X/Z)$ the set of classes in \cite[Theorem~18.2]{lyu-takumi-2026relativeminimalmodelprogram}.
    Then
    \[ 0 = [K_X+\Gamma] + [A] \in \mcA(X/Z)\subseteq \N^1(X/Z)_{\RR}. \]
    By \cite[Theorem~18.2(i)]{lyu-takumi-2026relativeminimalmodelprogram}, there exists a closed convex rational polytope
    $P\subseteq \N^1(X/Z)_{\RR}$ containing $0$ in its interior such that
    $P\cap \Nef(X/Z)$ is a closed convex rational polytope with nonempty interior.
    Hence $\Nef(X/Z)$ is rational polyhedral near $0$.
    As $\Nef(X/Z)$ is a cone, we deduce that $\Nef(X/Z)$ is a rational polyhedral cone.
    Since $\NE(X/Z)$ is the dual cone of $\Nef(X/Z)$ (cf. \cite[\S4.3]{lyu-takumi-2026relativeminimalmodelprogram}),
    we get that $\NE(X/Z)$ is a rational polyhedral cone.
\end{proof}


\subsection{Adjunction, complements, and normalisation}

We use the canonical bundle formula and the associated
comparison of boundary coefficients in
\cite[Lemma~2.38 and Theorem~2.39]{birkar2026steindegreefieldscharacteristic},
the relative complement theorem
\cite[Theorem~5.10]{birkar2026steindegreefieldscharacteristic},
and the estimate for degrees under normalisation in
\cite[Proposition~7.7]{birkar2026steindegreefieldscharacteristic}.
These results will be cited at the points where they are applied.


\subsection{\'Etale extensions}
The following lemma realises finite extensions of residue fields
by \'etale neighbourhoods while preserving the irreducibility
and normality of the total space after base change.

\begin{lemma}\label{et-field-extension}
    Let $X\to Z$ be a contraction between irreducible normal schemes of characteristic zero.
    Let $z\in Z$ be a scheme-theoretic point.
    Denote by $k(z)$ the residue field of $Z$ at $z$.
    Let $K/k(z)$ be a finite extension.
    Then there exist an integral scheme $U$ and a point $u\in U$ such that
    \begin{itemize}
        \item $U$ admits an \'etale morphism $U\to Z$,
        \item the image of $u\in U$ is $z$, 
        \item the residue field of $U$ at $u$ is isomorphic (over $k(z)$) to $K$, and
        \item $U\times_Z X$ is an irreducible normal scheme.
    \end{itemize}
\end{lemma}

\begin{proof}
    As $K/k(z)$ is a finite separable extension, it is simple \cite[Theorem~3, page~78]{Lorenz-alg-I}, 
    that is, there exists a monic separable irreducible polynomial $\bar{P}(T)\in k(z)[T]$ so that 
    \[ K\cong k(z)[T]/(\bar{P}(T)), \]
    where $T$ is an independent variable.
    Let $P(T)$ be a monic polynomial with coefficients in $\mcO_{Z,z}$
    so that $\bar{P}(T)$ lifts to $P(T)$.
    Then there is an affine open neighbourhood $V\subseteq Z$ of $z$ so that the coefficients of $P(T)$
    are regular functions on $V$.  As $\bar{P}(T)$ is separable, 
    the derivative $\bar{P}'(T)$ is invertible in $k(z)[T]/(\bar{P}(T))$.
    Hence up to shrinking $V$ near $z$, we can assume that $P'(T)$ is also invertible in $\mcO_Z(V)[T]/(P(T))$.  Set
    \[ U \coloneqq \spec \big(\mcO_Z(V)[T]/(P(T)) \big). \]
    By \cite[Proposition~18.38]{GW-agbook-2}, $U\to V$ is \'etale.
    Denote by $\mf p\subset \mcO_Z(V)$ the prime ideal corresponding to $z\in V$.  Then
    \[ \mf{q} \coloneqq \frac{\mf{p}[T] + (P(T))}{(P(T))} \]
    is a prime ideal of $\mcO_U(U)$.  Let $u\in U$ be the scheme-theoretic point
    corresponding to $\mf{q}$.  Then $u$ maps to $z\in V$, 
    and the residue field $k(u)$ of $U$ at $u$ is isomorphic to $K$.

    Since $U\to V$ is finite and \'etale, and since $V$ is normal and integral, 
    the irreducible components of $U$ are open and closed and are normal. 
    Replace $U$ by the irreducible component containing $u$. 
    Let $f_U\colon X_U\coloneqq U\times_Z X\to U$ be the base change. 
    Then flat base change gives
    \[ (f_U)_*\mcO_{X_U} = \mcO_U. \]
    Thus, $f_U$ has geometrically connected fibres. 
    Moreover, as $X_U\to X$ is \'etale, $X_U$ is normal. 
    Hence irreducible components of $X_U$ are open and closed. 
    Since $f_U$ is proper and surjective, at least one of the 
    finitely many irreducible components of $X_U$ dominates $U$, 
    and that irreducible component maps onto $U$. 
    Any additional irreducible component of $X_U$ produces a disconnected fibre of $f_U$. 
    Consequently, $X_U$ is irreducible.
\end{proof}


\subsection{Plt blow-ups}
The following lemma is a counterpart of \cite[Lemma~1]{Xu-fundamental-group} over arbitrary
fields of characteristic zero.
It produces a plt blow-up near the generic point of a prescribed
centre, allowing us to apply adjunction to the extracted divisor.

\begin{lemma}\label{plt-blowup}
    Let $(X,D)$ be an lc pair over a field $k$ of characteristic
    zero, and let $Z\subsetneq X$ be an integral closed subscheme
    of codimension $\ge 2$.  Assume that $(X,D)$ is klt
    along $Z$.
    Then up to shrinking $X$ around the generic point of $Z$,
    there exist an effective $\RR$-Cartier $\RR$-divisor
    $H$ on $X$, a birational contraction
    \[ \pi\colon Y\to X \]
    from a normal $k$-variety, and a $\QQ$-Cartier prime
    divisor $E$ on $Y$ with the following properties:
    \begin{itemize}
        \item $\Ex(\pi)= \Supp E$, $\pi(E) = Z$, and $-E$ is ample$/X$.
        \item The pair $(X,D+H)$ is lc, $E$ is its unique lc place over $X$, and $Z$ is its unique lc centre.
        \item The pair $(Y, D_Y+H_Y+E)$ is plt, and
              \[ K_Y+D_Y+H_Y+E = \pi^*(K_X+D+H), \]
               where 
              \[ D_Y \coloneqq \pi_*^{-1}D~~\mbox{and}~~ H_Y\coloneqq \pi_*^{-1}H. \]
        \item $-(K_Y+D_Y+E)$ is ample$/X$.
    \end{itemize}   
    Moreover, if $D$ has rational coefficients, then $H$ can
    be chosen to be a $\QQ$-Cartier $\QQ$-divisor.
\end{lemma}

\begin{proof}
    Shrinking $X$ around the generic point $\eta_Z$ of $Z$, we can
    assume that $X$ is affine and $(X,D)$ is klt.

    \medskip
    \emph{Step~1}.
    Let $\mc{I}_Z$ be the ideal sheaf of $Z$.  Take a 
    log resolution $f\colon W\to X$ of $(X,D)$ which also principalises
    $\mc{I}_Z$.  Write
    \[ K_W+D_W=f^*(K_X+D),~\mc{I}_Z\mcO_W=\mcO_W(-A),~\mbox{and}~ A=\sum_{i=1}^r a_i F_i, \]
    where $F_1,\dots,F_r$ are the exceptional prime divisors of $f$ and $a_i\in\ZZ^{\ge 0}$. 
    The support of $A$ is exceptional as $\codim(Z, X)\ge 2$.
    After replacing $X$ by a smaller affine neighbourhood of $\eta_Z$,
    every $F_i$ with $a_i>0$ has centre $Z$ on $X$. 
    The centre of an exceptional$/X$ prime divisor with $a_i = 0$ on $X$ may not be contained in $Z$.
    Set
    \[ \alpha_i \coloneqq \coeff_{F_i} D_W < 1~~\mbox{and}~~ I\coloneqq \Set{i\mid a_i>0}. \]
    It is evident that $I\neq\emptyset$.

    By \cite[1.12.1]{Kol_singularities_of_MMP},
    we can assume that there exists an effective exceptional$/X$ $\QQ$-divisor
    \[ G = \sum_{i=1}^r \beta_i F_i~~\mbox{with}~~ \beta_i\in\QQ^{>0} \]
    such that $-G$ is ample$/X$.
    By the openness of the ample cone \cite[Theorem~4.16]{lyu-takumi-2026relativeminimalmodelprogram}, we may perturb the $\beta_i$
    rationally so that the numbers $\beta_i/a_i$, for $i\in I$, are
    pairwise distinct, while keeping $\beta_i>0$ and $-G$ ample$/X$.
    Pick a sufficiently small $\epsilon\in \QQ^{>0}$ such that
    \[ \alpha_i + \epsilon\beta_i<1~\mbox{for every}~i, \]
    and such that the numbers
    \[ t_i\coloneqq \frac{1-\alpha_i-\epsilon\beta_i}{a_i},~i\in I, \]
    are pairwise distinct.  This is possible since the corresponding
    affine functions of $\epsilon$ have distinct slopes, so only
    finitely many values of $\epsilon$ are excluded.
    Let
    \[ t_0 \coloneqq \min_{i\in I} t_i>0. \]
    Up to reindexing, we can assume that this minimum is attained at $i=1$.  Then
    \[ \alpha_1 + t_0 a_1+\epsilon\beta_1=1~~\mbox{and}~~ \alpha_i+t_0a_i+\epsilon\beta_i<1~\mbox{for}~i\neq1. \]

    Pick an integer $m > \max \Set{1, t_0}$.
    As $X$ is affine, $\mc{I}_Z$ is generated by finitely many global sections.
    Let $L_1,\dots,L_m$ be the effective Cartier divisors defined
    by sufficiently general linear combinations of these generators.
    The moving parts of their pullbacks to $W$ belong to a
    base point free linear system.  Write
    \[ f^*L_j = A + \wt{L}_j,~ L\coloneqq \frac{1}{m}\sum_{j=1}^m L_j,~\mbox{and}~
    \wt L \coloneqq \frac{1}{m}\sum_{j=1}^m\wt{L}_j. \]
    Then by construction, we get
    \[ f^* L = A + \wt{L}. \]
    By Bertini's theorem \cite[Corollary~3.4.9]{joins-intersections}, we
    can assume $\wt{L}_j$ are reduced, they have no common irreducible components
    with each other or with 
    \begin{equation}\label{excep-divisor}
        (\Supp D_W) \cup (\bigcup_{1\le i\le r} F_i),
    \end{equation}
    and their union with the divisor \eqref{excep-divisor} is snc.
    Hence each irreducible component of $t_0\wt{L}$ has coefficient $t_0/m < 1$.    

    As $X$ is affine, by \cite[Corollaire~(4.6.6)]{EGA-II}, $-\epsilon G$
    is also globally ample on $W$.
    Pick a sufficiently large divisible $q\in \NN$ and a
    sufficiently general member $M$ of a finite dimensional very
    ample linear subsystem of $\abs{-q\epsilon G}$.  Set
    \[ H_0\coloneqq \frac{1}{q} M\sim_{\QQ}-\epsilon G. \]
    Let
    \[ \Sigma \coloneqq \Supp D_W \cup\bigcup_{i=1}^r F_i \cup\bigcup_{j=1}^m\Supp\wt{L}_j, \]
    which is an snc divisor by construction.  
    Again by Bertini's theorem, we can choose $M$ reduced, 
    with no common irreducible component with $\Sigma$, and transverse to all strata of $\Sigma$. 
    Then $\Sigma\cup\Supp M$ is snc.  In the boundary
    \[ D_W + t_0 f^*L + \epsilon G + H_0, \]
    the divisor $F_1$ has coefficient one, whereas every other irreducible component 
    has coefficient $<1$; indeed, the moving irreducible components have coefficients $t_0/m<1$ and $1/q<1$. 
    Hence the sub-pair
    \[ (W , D_W + t_0f^*L + \epsilon G + H_0) \]
    is sub-plt, with $F_1$ as its unique lc place over $X$.  Set
    \[ H_1 \coloneqq f_* H_0~~\mbox{and}~~ H\coloneqq t_0 L + H_1. \]
    As $H_0 + \epsilon G\sim_{\QQ}0$ and $G$ is $f$-exceptional, $H_1$ is effective and $\QQ$-Cartier. 
    Moreover, 
    \[ f^* H_1 = H_0 + \epsilon G. \]
    Consequently, $H$ is an effective $\RR$-Cartier $\RR$-divisor and
    \[ K_W+D_W+t_0f^*L+\epsilon G+H_0 = f^*(K_X+D+H). \]
    Thus, $(X,D+H)$ is lc, $F_1$ is its unique lc place, and $Z$ is its unique lc centre.  

    \medskip
    \emph{Step~2}.
    By \cite[Lemma~3.14]{birkar2026steindegreefieldscharacteristic},
    there exists a birational contraction $g\colon V\to X$
    from a $\QQ$-factorial $k$-variety $V$ such that
    the valuation of $F_1$ is a prime divisor $E_V$.
    Every other exceptional prime divisor has log
    discrepancy zero with respect to $(X,D+H)$, contradicting
    the uniqueness of $F_1$.  Thus, $E_V$ is the only exceptional
    prime divisor of $g$, and
    \[ K_V+\Lambda_V = g^*(K_X+D+H)~~\mbox{with}~~\Lambda_V=D_V+H_V+E_V, \]
    where $D_V$ and $H_V$ are the birational transforms of $D$ and $H$ to $V$.
    By construction in Step~1,
    it is evident that $(V,\Lambda_V)$ is plt.

    Pick $0<\delta<1$, and set 
    \[ B_V \coloneqq \Lambda_V - \delta E_V \ge 0. \]
    Then $(V,B_V)$ is klt, and
    \[ K_V+B_V \sim_{\RR} -\delta E_V/X. \]
    As $g$ is birational, $K_V+B_V$ is big over $X$.
    Then by \cite[Theorem~3.8]{birkar2026steindegreefieldscharacteristic},
    we can run a $(K_V+B_V)$-MMP over $X$ with scaling of a sufficiently
    general ample$/X$ divisor that terminates with a good minimal model
    $V'\to X$.  This is also a $(-E_V)$-MMP over $X$.
    By the negativity lemma, $E_V$ is not contracted by this MMP$/X$.
    Let $E_{V'}$ be the pushdown of $E_V$ to $V'$.
    Then $-E_{V'}$ is semiample$/X$. Let
    \[ \varphi\colon V'\to Y \]
    be the contraction over $X$ induced by $-E_{V'}$, and let
    \[ \pi\colon Y\to X \]
    be the induced birational contraction.
    Again by the negativity lemma, $E_{V'}$ is not contracted by $\varphi$.
    Hence $\varphi$ is a small contraction.  Let 
    \[ E \coloneqq \varphi_* E_{V'}, \]
    which is a $\QQ$-Cartier divisor on $Y$ such that $-E$ is ample$/X$,
    and $\pi(E) = Z$.
    (Note that $Y$ may not be $\QQ$-factorial.)
    Moreover, as $-E$ is ample$/X$, we get $\Ex(\pi) = \Supp E$.

    Denote by $D_Y, H_Y$ the birational transforms of $D, H$ to $Y$ respectively.
    Since $E$ is the only exceptional prime divisor of $\pi$,
    and since $E$ is an lc place of $(X,D+H)$,
    \[ K_Y+D_Y+H_Y+E = \pi^*(K_X+D+H). \]
    Then the uniqueness of the lc place implies that $(Y,D_Y+H_Y+E)$ is plt.
    Write
    \[ \pi^*H = H_Y + bE, \]
    where $b > 0$ is the order of $H$ along $E$.
    In particular, $H_Y$ is $\RR$-Cartier.  Then
    \[ -(K_Y+D_Y+E) \sim_{\RR} H_Y/X \sim_{\RR}-bE/X \]
    is ample$/X$.

    \medskip
    \emph{Step~3}.
    Finally, if $D$ has rational coefficients, then all $\alpha_i$'s in Step~1
    are rational.  Since $\epsilon$, the $\beta_i$, and the $a_i$ in Step~1
    are rational, we deduce that $t_0$ is also rational.  Thus, the divisor $H$
    is a $\QQ$-Cartier $\QQ$-divisor.
\end{proof}


The following lemma provides the Fano type contraction needed to apply induction after adjunction.

\begin{lemma}\label{plt-induction}
    Let $k$ be a field of characteristic zero, and let
    $f\colon Y\to Z$ be a birational contraction of Fano type over $k$.
    Let $(Y, B_Y)$ and $(Z, B_Z)$ be pairs satisfying
    \[ K_Y+B_Y = f^*(K_Z+B_Z). \]
    Assume that $(Y,B_Y)$ is plt and that $E\coloneqq \rddown{B_Y}$ is a
    prime divisor.  Let $T\coloneqq f(E)$ with its reduced structure, and let $k_E$
    be the algebraic closure of $k$ in $K(T)$.  Then $E\to T$ is a
    contraction of Fano type between $k_E$-varieties.
\end{lemma}

\begin{proof}
    \emph{Step~1}.
    Denote by $g$ the induced morphism $E\to T$.
    We first show that $g$ is a contraction.
    Let $\bar k$ be the algebraic closure of $k$, and let
    $\bar f\colon \bar{Y}\to \bar{Z}$ be the induced morphism, where 
    $\bar{Y}\coloneqq Y\otimes_k \bar{k}$ and $\bar{Z}\coloneqq Z\otimes_k \bar{k}$.
    Then $\bar{f}$ is a birational contraction of Fano type between $\bar{k}$-varieties;
    see \cite[Lemma~2.47]{birkar2026steindegreefieldscharacteristic}.
    Let $B_{\bar{Y}}$ be the flat pullback of $B_Y$ to $\bar{Y}$.
    By \cite[\S 2.35(5)]{birkar2026steindegreefieldscharacteristic},
    the non-klt locus of $(\bar{Y}, B_{\bar{Y}})$ is $\bar{E}\coloneqq E\otimes_k \bar{k}$.
    Set $\bar{T}\coloneqq T\otimes_k \bar{k}$.
    Then by \cite[Theorem~1.3(2)]{birkar-connected-JDG}, 
    the fibres of $\bar{E}\to \bar{T}$ are connected.
    Hence $E\to T$ has geometrically connected fibres.
    Let 
    \[ E\to V\xrightarrow{\nu} T \]
    be the Stein factorisation of $E\to T$.
    Then $\nu$ is finite, surjective, and radicial; see \cite[Chap.~I, (3.5.5)]{EGA-I}.
    As $\chara k = 0$, the residue field extensions of $V\to T$ are trivial.
    By \cite[(11.10.2)]{Kol-general-type}, $T$ is seminormal.
    Thus, by \cite[Definition~10.74]{Kol-general-type}, $\nu$ is an isomorphism.
    Then we get $g_* \mcO_E = \mcO_T$, i.e. $g$ is a contraction.
    Moreover, as $(Y, B_Y)$ is plt, $E$ is normal \cite[Theorem~4.16]{Kol_singularities_of_MMP}.
    Hence $T$ is also normal.

    By \cite[Proposition~(4.5.9)]{EGA-IV-2}, $K(T)$ is algebraically closed in $K(E)$.
    Then $k_E$ is also the algebraic closure of $k$ in $K(E)$.
    Hence $g\colon E\to T$ is a contraction between $k_E$-varieties.

    \medskip
    \emph{Step~2}.
    By \cite[\S 2.46]{birkar2026steindegreefieldscharacteristic}, 
    there exists an effective $\QQ$-divisor $G$ on $Y$ such that $-(K_Y+G)$ is ample$/Z$.
    Pick an effective Cartier divisor $H$ on $Z$ containing $T$.  Set
    \[ a\coloneqq \coeff_E G < 1~~\mbox{and}~~ m\coloneqq \coeff_E f^*H > 0. \]
    Moreover, let
    \[ C\coloneqq G + \frac{1-a}{m}f^*H. \]
    Then $C \ge 0$ and $\coeff_E C = 1$.  For a sufficiently small $0< \epsilon \ll 1$, write
    \[ \Delta_\epsilon\coloneqq (1-\epsilon)B_Y + \epsilon C. \]
    On the sufficiently high log resolution $\pi\colon W\to Y$ of $(Y, B_Y + C)$, the coefficient of 
    the birational transform of $E$ in $\pi^*(K_Y + \Delta_{\epsilon}) - K_W$ is one.  
    As $(Y, B_Y)$ is plt, all other coefficients of  
    of $\pi^*(K_Y + \Delta_{\epsilon}) - K_W$ is $<1$.
    Hence $(Y, \Delta_{\epsilon})$ is plt and $\rddown{\Delta_{\epsilon}} = E$.
    Moreover,
    \[ -(K_Y + \Delta_{\epsilon})\sim_{\RR} -\epsilon (K_Y + G)/Z \]
    is ample$/Z$.
    By divisorial adjunction \cite[\S 11.1]{Kol-general-type}, we can write
    \[ K_E + \Gamma = (K_Y+\Delta_{\epsilon})|_E. \]
    Moreover, by \cite[Theorem~11.17]{Kol-general-type}, $(E, \Gamma)$ is klt.
    As $-(K_Y + \Delta_{\epsilon})$ is ample$/Z$, we deduce that 
    $-(K_E + \Gamma)$ is ample$/T$.  This shows that 
    $g\colon E\to T$ is a contraction of Fano type between $k_E$-varieties.
\end{proof}


\section{\bf First proof of the main result}\label{s-first-pf}

In this section, we prove Theorem~\ref{dim-d-t} using the
result \cite[Theorem~1.3]{Hogadi-Xu-09}, together with plt blow-ups and adjunction.

\begin{theorem}[={Theorem~\ref{dim-d-t}}]\label{bir-pf}
    Let $d\in \NN$ and let $t\in (0,1]$.
    There exists an $N_b(d,t) \in \NN$ depending only on $d,t$ such that
    the following holds.

    Let $k$ be a field of characteristic zero.
	Let $(X, B + \mb{M})\to Z$ be a generalised log Calabi-Yau fibration over $k$,
	where $f\colon X\to Z$ is a birational contraction of Fano type between $k$-varieties with $\dim X = d$.
    Let $S$ be an irreducible component of $B$ 
	whose coefficient in $B$ is $\ge t$.  Then 
    \[ \sdeg (S^{\nor}/T) \le N_b(d,t), \]
    where $T$ is the image of $S$ in $Z$, and $S^{\nor}$ is the normalisation of $S$.
\end{theorem}

\begin{proof}
    By \cite[Lemmas~2.50 and~3.14]{birkar2026steindegreefieldscharacteristic}, 
    we can assume that $(X,B+\mb{M})$ is a $\QQ$-factorial
    g-dlt g-pair.  Denote by $\eta_T$ the generic point of $T$.
    If $S$ is not exceptional over $Z$, then $S\to T$ is birational and the required Stein degree is one. 
    We therefore assume that $\codim(T, Z)\ge 2$; in particular, $d\geq2$. 
    Throughout the proof, $t$ denotes the coefficient lower bound in the statement of the theorem.

    \medskip
    \emph{Step~1}.
    In this step, we reduce to the case of ordinary pairs with rational coefficients.  Let
    \[ b\coloneqq \coeff_S B \in \RR^{>0}. \]
    By \cite[Lemma~3.9]{birkar2026steindegreefieldscharacteristic}, we can run a $(-K_X - bS)$-MMP $\phi\colon X\bir Y$
    over $Z$ which terminates with a good minimal model $Y\to Z$ of $-(K_X+bS)$.
    Denote by
    \[ X = X_0\xdashrightarrow{\phi_0} X_1\xdashrightarrow{\phi_1} \cdots \bir X_i\xdashrightarrow{\phi_i} \cdots \]
    the sequence of birational contractions in this MMP$/Z$.
    
    We claim that $S$ is not contracted by $\phi$.
    To this end, assume that $S$ is contracted by $\phi_i\colon X_i\bir X_{i+1}$.
    Denote by $S_i$ the pushdown of $S$ to $X_i$.  Then 
    $\phi_i$ is a divisorial contraction whose exceptional locus is equal to the support of $S_i$.
    Let $p\colon W\to X_i$ be a sufficiently high log resolution so that $\mb{M}$
    descends to a nef$/Z$ $\RR$-divisor $\mb{M}_W$ on $W$, and let $q\colon W\to X_{i+1}$
    be the composite birational contraction $\phi_i\circ p$.
    By the negativity lemma \cite[Lemma~2.24]{birkar2026steindegreefieldscharacteristic}, we can write 
    \[ p^* \mb{M}_{X_i} = \mb{M}_W + E, \]
    where $E\ge 0$ is exceptional$/X_i$.
    Pick a closed point of $S_i$ outside $\Supp (B_i-bS_i)\cup p(\Ex (p))$ 
    and a curve $C$ contracted by $\phi_i$ through this closed point.
    As $\phi_i$ is a divisorial contraction with exceptional divisor $S_i$,
    such curves cover $S_i$.
    Let $C_W$ be the birational transform of $C$ on $W$.
    Then $C_W$ is not contained in $\Supp E$.  Thus, the projection formula gives
    \[ \mb{M}_{X_i}\cdot C =(\mb{M}_W+E)\cdot C_W\ge 0. \]
    As $B_i-bS_i$ is effective and $C$ is not contained in $\Supp (B_i - bS_i)$, this contradicts
    \[ (B_i - bS_i + \mb{M}_{X_i})\cdot C = (-K_{X_i} - b S_i)\cdot C < 0. \]
    Thus, we can conclude that $S$ is not contracted by $\phi$.
    Denote by $S_Y$ the pushdown of $S$ to $Y$.

    Assume that $b$ is an irrational number.  
    Let $R$ be the extremal ray of $\NE(X_i/Z)$ that is contracted by $X_i\bir X_{i+1}$.
    Then 
    \[ -(K_{X_i} + bS_i)\cdot (R\setminus\Set{0}) < 0. \]
    Thus, for any $b'\in (0,b)\cap \QQ$ that is sufficiently close to $b$, we also have
    \[ -(K_{X_i} + b'S_i)\cdot (R\setminus\Set{0}) < 0. \]
    For the resulting model $Y\to Z$ of $(-K_X- bS)$-MMP, by Lemma~\ref{fano-type-polyhedral}, 
    $\NE(Y/Z)$ is a rational polyhedral cone.  Denote by $R_1, \dots, R_l$
    the extremal rays of $\NE(Y/Z)$, which are defined over $\QQ$.  If
    \[ -(K_Y + bS_Y)\cdot (R_j\setminus \Set{0}) > 0, \]
    then for any $b'\in (0,b)\cap \QQ$ sufficiently close to $b$, the inequality also holds
    for $-(K_Y + b'S_Y)$.  If
    \[ -(K_Y + bS_Y)\cdot R_j = 0, \]
    then as $b$ is an irrational number, we see that
    \[ K_Y \cdot R_j = 0~~\mbox{and}~~ S_Y\cdot R_j = 0. \]
    Thus, for any $b'\in (0,b)$, we also have
    \[ -(K_Y+b'S_Y) \cdot R_j = 0. \]
    Hence for any $b'\in (0,b)\cap \QQ$ sufficiently close to $b$, $-(K_Y + b'S_Y)$
    is also nef$/Z$.  Therefore, $\phi\colon X\bir Y$ is also a $(-K_X-b'S)$-MMP over $Z$, 
    which terminates with a good minimal model $Y\to Z$ of $(-K_X-b'S)$; see \cite[Corollary~2.55]{birkar2026steindegreefieldscharacteristic}.
    Moreover, by \cite[Lemma~2.32]{birkar2026steindegreefieldscharacteristic}, 
    $(Y, b' S_Y)$ is an lc pair over $k$.

    If $b$ is rational, we set $b'=b$; if $b$ is irrational, we pick a $b'\in (0,b)\cap \QQ$
    sufficiently close to $b$ so that $b'> t/2$.  
    Pick a closed point $z\in T$.
    By \cite[Theorem~5.10]{birkar2026steindegreefieldscharacteristic},
    up to shrinking $Z$ to an open neighbourhood of $z$ containing $\eta_T$,
    there exist an $n\in \NN$ depending only on $d,b'$ and a monotonic $n$-complement $K_Y + C_Y$ of $K_Y + b'S_Y$.
    In particular, we get 
    \[ K_Y + C_Y\sim_{\QQ}0/Z~~\mbox{and}~~ \coeff_{S_Y} C_Y \ge \frac{t}{2}. \]
    Thus, replacing $t$ with $t/2$ and $(X, B+\mb{M})\to Z$ with $(Y, C_Y)\to Z$,
    we can assume that $\mb{M} = 0$, and $B$ has rational coefficients.

    \medskip
    \emph{Step~2}.
    Since $X\to Z$ is of Fano type, and since $X$ is $\QQ$-factorial, 
    by \cite[Lemma~3.9]{birkar2026steindegreefieldscharacteristic}, 
    we can run a $(-S)$-MMP over $Z$ that terminates with a good minimal model $X'\to Z$ of $-S$.
    Denote by $S'$ the pushdown of $S$ to $X'$.
    Then $-S'$ is semiample$/Z$.  
    Let $X'\to X''$ be the contraction over $Z$ induced by $-S'$,
    and let $S''$ be the pushdown of $S'$ to $X''$.
    Then $-S''$ is ample$/Z$.
    Hence up to shrinking $Z$ around $\eta_T$, $\Supp S''$ is equal to
    the inverse image of $T$ in $X''$.
    Thus, replacing $X$ with $X''$ and shrinking $Z$ around $\eta_T$, 
    we can assume that $S$ is the only irreducible component of $f^{-1}(T)$.
    Note that after this process, $X$ may not be $\QQ$-factorial;
    however, $S$ is $\QQ$-Cartier.

    As $X\to Z$ is of Fano type, there exists a $\QQ$-divisor $\Delta$ on $X$ such that
    $(X, \Delta)$ is klt and $K_X+\Delta\sim_{\QQ}0/Z$.  Set $\Delta_Z \coloneqq f_* \Delta$.  Then
    \[ K_X+\Delta = f^*(K_Z + \Delta_Z). \]
    Thus, $(Z, \Delta_Z)$ is klt.  
    Let $k(\eta_T)$ be the residue field of $Z$ at $\eta_T$, and let
    $X_{\eta_T}$ be the scheme-theoretic fibre of $X\to Z$ over $\eta_T$.
    By \cite[Theorem~1.3]{Hogadi-Xu-09}, there exists a geometrically integral subscheme $G$
    of $f^{-1}(\eta_T)$, which is defined over $k(\eta_T)$.
    Let $K/k(\eta_T)$ be a finite Galois extension such that
    each irreducible component of $S_{\eta_T}\otimes_{k(\eta_T)} K$ is geometrically irreducible.
    Denote by
    \[ S_{K, 1}, \dots, S_{K, r} \]
    the irreducible components of $S_{\eta_T}\otimes_{k(\eta_T)} K$.
    By Lemma~\ref{et-field-extension}, there exists an irreducible normal scheme $U$
    admitting an \'etale morphism $U\to Z$ such that a scheme-theoretic point $u\in U$
    maps to $\eta_T$ and that $k(u)\cong K$.  
    Shrinking $U$ around $u$, we can assume that $u$ is the inverse image of $\eta_T$ in $U$.
    Note that $G\otimes_{k(\eta_T)} K$ is also irreducible.
    As $\Gal(K/k(\eta_T))$ acts transitively, we have
    \[ G\otimes_{k(\eta_T)} K \subseteq S_{K, i} \]
    for each $1\le i\le r$.  If $\dim S_{\eta_T} = \dim G$,
    then $\sdeg (S^{\nor}/T) = 1$.  Hence in the remainder of the proof,
    we assume that $\dim G < \dim S_{\eta_T}$.

    \medskip
    \emph{Step~3}.
    Let $\Delta$ be the $\QQ$-divisor in Step~2, and let
    \[ \Gamma\coloneqq \frac{1}{2}\Delta + \frac{1}{2}B. \]
    Then $(X, \Gamma)$ is also klt, and $K_X+\Gamma \sim_{\QQ} 0/Z$.
    Moreover, it is evident that $\coeff_S \Gamma \ge t/2$.
    
    Let $R$ be the closure of $G$ in $X$ and $\eta_R$ the generic point of $R$.
    By Lemma~\ref{plt-blowup}, up to shrinking $X$ around $\eta_R$,
    there exist a $\QQ$-Cartier $\QQ$-divisor $H$ on $X$ whose support contains $R$
    and a birational contraction $\pi\colon Y\to X$ from a normal $k$-variety such that
    \begin{itemize}
        \item $Y$ has a prime divisor $E$ with $\centre_X E = R$,
        \item $-E$ is ample$/X$, 
        \item $(Y, \pi_*^{-1}\Gamma + E + H_Y)$ is plt,
        \item $R$ is the unique lc centre of $(X, \Gamma+H)$, and
        \item $E$ is the unique lc place of $(X, \Gamma + H)$,
    \end{itemize}
    where $H_Y$ is the birational transform of $H$ to $Y$.
    Denote by $S_Y$ the birational transform of $S$ to $Y$.

    As $(Y, \pi_*^{-1} \Gamma + E + H_Y)$ is plt, $E$ is normal \cite[Theorem~4.16]{Kol_singularities_of_MMP}.
    Moreover, as $\pi$ is a contraction, the generic fibre of $E\to R$
    is geometrically normal and geometrically connected, whence geometrically integral.
    Thus, $E\to R$ is a contraction.
    Since $-E$ is ample$/X$, $S_Y$ is ample$/X$.
    Hence the reduction of $S_Y\cap E$ is a divisor on $E$;
    indeed, the restriction of $S_Y$ to the positive-dimensional generic
    fibre of $E\to R$ is an effective ample divisor, hence nonzero.
    Let $V$ be an irreducible component of $S_Y \cap E$ that dominates $R$.
    By adjunction \cite[\S 4.1]{Kol_singularities_of_MMP}, we can write
    \[ K_E + \Theta_E = (K_Y + \pi_*^{-1} \Gamma + E + H_Y)|_E, \]
    where $(E, \Theta_E)$ is klt, and $\coeff_V \Theta_E \ge t/2$.
    Moreover, as $K_Y + \pi_*^{-1} \Gamma + E + H_Y$ is equal to $\pi^*(K_X+\Gamma + H)$, we get
    \[ K_E+\Theta_E \sim_{\QQ}0/R. \]
    Note that $E\to T$ is dominant.  Denote by $E_{\eta_T}$ the generic fibre of $E\to T$.
    Then the generic fibre of $E_{\eta_T}\to G$ is geometrically integral over $K(R) = K(G)$.
    As $G$ is geometrically irreducible over $k(\eta_T)$,
    we can conclude that $E_{\eta_T}$ is also geometrically irreducible over $k(\eta_T)$.
    Hence $E_{\eta_T}\to G$ is a contraction between $k(\eta_T)$-varieties.
    Denote by $V_{\eta_T}$ the generic fibre of $V\to T$, which is irreducible.
    By construction as above, 
    \[ (E_{\eta_T}, (\Theta_E)|_{E_{\eta_T}})\to G \]
    is a log Calabi-Yau fibration over $k(\eta_T)$ with 
    $(E_{\eta_T}, (\Theta_E)|_{E_{\eta_T}})$ klt, and 
    \[ \coeff_{V_{\eta_T}} (\Theta_E)|_{E_{\eta_T}} \ge t/2. \]
    Then by \cite[Theorem~1.1]{birkar2026steindegreefieldscharacteristic},
    \[ \ell \coloneqq \sdeg(V^{\nor}/R) = \sdeg (V_{\eta_T}^{\nor}/G) \]
    is bounded from above depending only on $d,t$.

    \medskip
    \emph{Step~4}.
    As $U\times_Z Y \to Y$ is \'etale, $U\times_Z Y$ is normal.
    Moreover, since $\pi$ is a birational contraction, $U\times_Z Y\to  U\times_Z X$
    is also a birational contraction.  Hence $Y_U\coloneqq U\times_Z Y$ is an irreducible normal scheme. 
    We include the following diagram for convenience.
    \[\xymatrix{
     U\times_Z Y\ar[r]\ar[d] & Y\ar[d]^{\pi} \\
     U\times_Z X\ar[r]\ar[d] & X\ar[d]^f \\
      U\ar[r] & Z
    }\]
    Set
    \[ E_U\coloneqq U\times_Z E ~~\mbox{and}~~ R_U \coloneqq U\times_Z R = \ol{\Set{u}}\times_T R. \]
    Up to shrinking $U$ around $u$, we can assume that $R_U$ is irreducible,
    $E_U$ is an irreducible normal scheme, and $E_U\to R_U$ is a contraction.

    By construction in Step~2, each $S_{K, i}$ is a geometrically irreducible scheme over $k(u)\cong K$,
    and the generic point of $S_{K, i}$ is a codimension one point of $X_U \coloneqq U\times_Z X$.
    Denote by $S_{U, i}$ the closure of $S_{K, i}$, and by $S_{U, i}^Y$ the
    birational transform of $S_{U, i}$ to $Y_U$.  Then 
    \[ R_U \subset S_{U,i}~\mbox{for each}~1\le i\le r. \]
    Let 
    \[ V_{K,1}, \dots, V_{K,s} \]
    be the irreducible components of $V_{\eta_T}\otimes_{k(\eta_T)} K$.
    It is evident that $s\le \ell$.
    Up to shrinking $U$ around $u$, we can assume that
    $V\times_R R_U$ also has $s$ irreducible components.  Let
    \[ V_{U, 1}, \dots, V_{U, s} \]
    be the irreducible components of $V\times_R R_U$ such that each $V_{U, j}$
    is the closure of $V_{K, j}$.    
    
    We claim that each $S_{U, i}^Y$ contains at least one $V_{U, j}$ for some $1\le j\le s$.
    Let $S_{Y, \eta_T}$ be the generic fibre of $S_Y \to T$, and let 
    \[ S_{K,1}^Y, \dots, S_{K, m}^Y \]
    be the irreducible components of $S_{Y, \eta_T}\otimes_{k(\eta_T)} K$.
    As $Y\to X$ is an isomorphism over the generic point of $S$, we get $m = r$,
    and up to reindexing, we can assume that $S_{K,i}^Y\to S_{K, i}$ is birational for each $1\le i\le r$.
    Since $\Gal(K/k(\eta_T))$ acts transitively on $\Set{S_{K,1}, \dots, S_{K,r}}$,
    the group $\Gal(K/k(\eta_T))$ also acts transitively on $\Set{S_{K,1}^Y, \dots, S_{K,r}^Y}$.
    Moreover, as $V_{\eta_T}\subset S_{Y, \eta_T}$, we see that for every $1\le j' \le s$, $V_{K, j'}$
    is contained in $S_{K, i'}^Y$ for some $1\le i' \le r$.
    Pick an arbitrary $S_{K, i}^Y$.  Let $\sigma\in \Gal(K/k(\eta_T))$ such that
    $\sigma(S_{K, i'}^Y) = S_{K, i}^Y$.  Then $\sigma(V_{K, j'})$ is contained in $S_{K, i}^Y$.
    As $\Gal(K/k(\eta_T))$ acts transitively on $\Set{V_{K,1}, \dots, V_{K, s}}$,
    the closure of $\sigma(V_{K, j'})$ is $V_{U, j}$ for some $1\le j\le s$.
    Then $V_{U, j}$ is contained in $S_{U, i}^Y$ as desired.   
    
    We claim that for every $V_{U, j}$, the number of divisors $S_{U, i}^Y$ containing $V_{U, j}$
    is bounded from above by some $N\in \NN$ depending only on $t$.
    Let $\Gamma_U$ be the flat pullback of $\pi^{-1}_* \Gamma$ to $Y_U$, and let 
    $E_U$ be the flat pullback of $E$ to $Y_U$.  As $E_U$ is irreducible,
    by \cite[\S 2.34]{birkar2026steindegreefieldscharacteristic}, $(Y_U, \Gamma_U + E_U)$ is plt.
    By adjunction \cite[Lemma~4.8]{Kol_singularities_of_MMP}, we can write
    \[ K_{E_U} + D_U = (K_{Y_U} + \Gamma_U + E_U)|_{E_U}, \]
    where $(E_U, D_U)$ is klt.  
    Denote by $I_j\subseteq \Set{1, \dots, r}$ the set of indices such that for each $i\in I_j$,
    we have $V_{U, j}\subset S_{U, i}^Y$.
    By \cite[Theorem~3.36]{Kol_singularities_of_MMP}, 
    \[ \dfrac{m-1 + (t/2)\sum_{i\in I_j} n_i}{m} \le \coeff_{V_{U, j}} D_U <1, \]
    where $m, n_i$ are positive integers.  Hence
    \[ (t/2)\abs{I_j} < 1, \]
    where $\abs{I_j}$ denotes the cardinality of $I_j$.  This proves the claim.

    By the two claims above, we can conclude that
    \[ \sdeg(S^{\nor}/T) = r \le s N \le \ell N, \]
    which is bounded from above depending only on $d, t$.
\end{proof}


\section{\bf Second proof of the main result}\label{s-second-pf}

In this section, we prove Theorem~\ref{thm-vertical-introduction}
and give a second proof of Theorem~\ref{dim-d-t} by simultaneous induction on dimension.


\begin{lemma}\label{vertical-to-birational}
    Let $d\in \NN$.
    Assume Theorem~\ref{thm-vertical-introduction} holds for $d-1$ and
    each $t\in (0,1]$.
    Then Theorem~\ref{dim-d-t} holds for $d$ and every $t\in (0,1]$.
\end{lemma}

\begin{proof}
    Let $k$ be a field of characteristic zero.
	Let $(X, B + \mb{M})\to Z$ be a generalised log Calabi-Yau fibration over $k$,
	where $f\colon X\to Z$ is a birational contraction of Fano type between $k$-varieties, and $\dim X = d$.
    Let $S$ be an irreducible component of $B$ 
	whose coefficient in $B$ is $\ge t$.  Denote by $T$ the image of $S$ in $Z$.
    By Step~1 in the proof of Theorem~\ref{bir-pf},
    we can assume that $\mb{M} =0$ and $B$ has rational coefficients.
    If $S$ is not exceptional$/Z$, then $S\to T$ is birational,
    so $\sdeg (S^{\nor}/T) = 1$.  Thus, we assume that $\codim (T,Z)\ge 2$.
    Moreover, if $d=1$, then $X\to Z$ is an isomorphism of nonsingular curves.
    Hence we also assume that $d\ge 2$.
    
    \medskip
    \emph{Step~1}.
    Set $B_Z \coloneqq f_* B$.  Since $f$ is birational, and since $K_X+B\sim_{\QQ}0/Z$, we get
    \[ K_X + B = f^*(K_Z + B_Z). \]
    Since $X$ is of Fano type over $Z$, there exists an effective
    $\QQ$-divisor $\Delta$ such that $(X, \Delta)$ is klt and
    $K_X+\Delta\sim_{\QQ}0/Z$.  Replacing $B$ by
    $\frac{1}{2}(B + \Delta)$, we can assume that $(X,B)$ is klt 
    up to replacing the current coefficient lower bound
    $t$ by $t/2$.  Then $(Z,B_Z)$ is klt.
    By Lemma~\ref{plt-blowup}, there exist a $\QQ$-Cartier $\QQ$-divisor $H_Z$ on $Z$ 
    and a birational contraction 
    \[ \pi\colon Y\to Z \] 
    from a normal $k$-variety $Y$ such that
    \begin{itemize}
        \item $Y$ has a prime divisor $E$ whose centre on $Z$ is equal to $T$,
        \item $-(K_Y + \pi^{-1}_*B_Z + E)$ and $-E$ are ample$/Z$,
        \item $T$ is the unique lc centre of $(Z, B_Z + H_Z)$, 
        \item $(Y, \pi^{-1}_* B_Z + E)$ is plt, and
        \item $E$ is the unique lc place of $(Z, B_Z + H_Z)$ on $Y$.
    \end{itemize}
    By \cite[Lemma~2.50]{birkar2026steindegreefieldscharacteristic}, $Y\to Z$ is of Fano type.
    It is evident that 
    \[ a(S, Z, B_Z + H_Z) \le a(S, Z, B_Z) \le 1-t. \]
    By Lemma~\ref{plt-induction}, $E\to T$ is a contraction, and $E$ is normal.
    If $S$ and $E$ correspond to the same geometric valuation, then
    $\sdeg (S^{\nor}/T) = \sdeg (E/T) = 1$.
    Hence we assume that $S$ and $E$ define distinct geometric valuations.
    
    Denote by $C$ the centre of $S$ on $Y$.
    As $E$ is the unique exceptional$/Z$ divisor of $\pi$,
    we get $C\subsetneq E$.
    By \cite[Lemma~3.14]{birkar2026steindegreefieldscharacteristic}, there exists a birational contraction
    $\phi\colon Y'\to Y$ extracting only the geometric valuation of $S$, where $Y'$ is $\QQ$-factorial.
    Denote by $S'$ the centre of $S$ on $Y'$, which is the only $\phi$-exceptional divisor on $Y'$.  Write
    \[ K_{Y'} + \Delta_{Y'} = (\pi\circ\phi)^*(K_Z + B_Z + H_Z). \]
    Let $E'$ be the birational transform of $E$ to $Y'$.   
    Then $(Y', \Delta_{Y'})$ is plt with $\rddown{\Delta_{Y'}} = E'$.
    Moreover, we have
    \[ \coeff_{S'} \Delta_{Y'} = 1 - a(S, Z, B_Z + H_Z) \ge t. \]
    By \cite[Lemma~2.50]{birkar2026steindegreefieldscharacteristic}, $Y'\to Z$ is of Fano type.
    
    We claim there exists an irreducible component of $S'\cap E'$
    that dominates $C$.  Indeed, as $-E$ is ample$/Z$, $E$ is $\QQ$-Cartier, so we can write
    \[ \phi^* E = E' + \lambda S' \]
    for some $\lambda\in \QQ^{>0}$.  Then by \cite[(21.4.7)]{EGA-IV-4}, we get
    \begin{equation}\label{set-theorectic-pullback}
        \phi^{-1}(E) = S'\cup E'
    \end{equation}
    set-theoretically. 
    Denote by $\eta_C$ the generic point of $C$.
    Since $E'\to E$ and $S'\to C$ are projective and surjective, 
    $\phi^{-1}(\eta_C)$ intersects both $E'$ and $S'$.
    By \eqref{set-theorectic-pullback}, a scheme-theoretic point of $\phi^{-1}(\eta_C)$
    is contained in either $S'$ or $E'$.  Thus, if $S'\cap E' = \emptyset$,
    then $S'\cap \phi^{-1}(\eta_C)$ and $E'\cap \phi^{-1}(\eta_C)$ give at least two 
    connected components of $\phi^{-1}(\eta_C)$, contradicting that $\phi$ is a contraction.
    Hence $S'\cap E'$ is nonempty, and $S'\cap E'$ intersects $\phi^{-1}(\eta_C)$.  
    As $S'$ is $\QQ$-Cartier, each irreducible component
    of $S'\cap E'$ is a prime divisor on $E'$.
    Pick an irreducible component $P$ of $S'\cap E'$ that intersects $\phi^{-1}(\eta_C)$.
    Then the image of $P$ in $Z$ is equal to $T$.

    \medskip
    \emph{Step~2}.
    By divisorial adjunction \cite[Definition~11.14]{Kol-general-type}, we can write
    \[ K_{E'} + \Delta_{E'} = (K_{Y'} + \Delta_{Y'})|_{E'}. \]
    By \cite[Theorem~11.17]{Kol-general-type}, $(E', \Delta_{E'})$ is klt.
    Moreover, by \cite[Lemma~11.16]{Kol-general-type}, we get
    \[ \coeff_P \Delta_{E'} \ge \coeff_{S'} \Delta_{Y'} \ge t. \]
    By Lemma~\ref{plt-induction}, $E'\to T$ is a contraction of Fano type between
    $k_{E'}$-varieties, where $k_{E'}$ denotes the algebraic closure of $k$ in $K(E')$.
    By construction, we get
    \[ K_{E'} + \Delta_{E'}\sim_{\QQ}0/T. \]
    Since $\dim E' = \dim Y' - 1 = d-1$, the hypothesis gives Theorem~\ref{thm-vertical-introduction} 
    over the field $k_{E'}$ in dimension $d-1$.
    Hence
    \[ \sdeg (P^{\nor}/T) \le N_{\mbb V}(d-1,t), \]
    which depends only on $d, t$.

    Let $(S')^{\nor}\to S'$ be the normalisation morphism,
    and let $Q$ be a prime divisor on $(S')^{\nor}$ mapping onto $P$.
    As $\coeff_{S'} \Delta_{Y'} \ge t$,
    by \cite[Proposition~7.7]{birkar2026steindegreefieldscharacteristic},
    the degree of $Q\to P$ is bounded from above depending only on $t$.
    As $Q^{\nor}\to T$ is surjective, by \cite[Lemmas~2.9, 7.1]{birkar2026steindegreefieldscharacteristic}, we get
    \[ \sdeg (S^{\nor}/T) = \sdeg ((S')^{\nor}/T) \le \sdeg (Q^{\nor}/T) 
    \le \deg (Q/P)\cdot \sdeg (P^{\nor}/T), \]
    which is bounded from above depending only on $d,t$.
    This shows that Theorem~\ref{dim-d-t} holds for $d$ and every $t\in (0,1]$.
\end{proof}


\begin{lemma}\label{birational-to-vertical}
    Let $d\in \NN$.  Assume Theorem~\ref{dim-d-t} holds for $d$
    and each $t\in (0,1]$.  Then Theorem~\ref{thm-vertical-introduction} holds for $d$
    and each $t\in (0,1]$.
\end{lemma}

\begin{proof}
    Let $k$ be a field of characteristic zero.  Let 
    \[ f\colon (X, B + \mb{M})\to Z \] 
    be a Fano type generalised log Calabi-Yau fibration over $k$ with $\dim X = d$,
    and let $S$ be an irreducible component of $B$ with $\coeff_S B \ge t$.
    Denote by $T$ the image of $S$ in $Z$.
    If $\codim(T, Z) \le 1$, then we can conclude by 
    \cite[Theorems~1.1, 1.2]{birkar2026steindegreefieldscharacteristic}.
    Thus, we can assume that $\codim (T, Z) \ge 2$.

    By \cite[Lemmas~2.50, 3.14]{birkar2026steindegreefieldscharacteristic},
    we can assume that $(X, B + \mb{M})$ is a $\QQ$-factorial
    g-dlt g-pair.  Then by \cite[Lemma~3.9]{birkar2026steindegreefieldscharacteristic},
    we can run a $(-S)$-MMP $\phi\colon X\bir X'$ over $Z$.
    As $-S$ is pseudo-effective$/Z$, the resulting model $X'\to Z$
    is a good minimal model of $-S$ over $Z$.
    By the negativity lemma \cite[Lemma~2.24]{birkar2026steindegreefieldscharacteristic},
    $S$ is not contracted by $\phi$.  Let $S'\coloneqq \phi_* S$.  
    Then $-S'$ is semiample$/Z$.  
    Let $g\colon X'\to Y'$ be the contraction over $Z$ induced by $-S'$.
    As $S'$ is vertical$/Z$, we deduce that $Y'\to Z$ is birational.
    Denote by $T'$ the image of $S'$ in $Y'$.  
    Then up to shrinking $Z$ around the generic point of $T$, $-T'$ is ample$/Z$.  

    Let $B'\coloneqq \phi_* B$.  Then we have 
    $K_{X'} + B' + \mb{M}_{X'}\sim_{\RR}0/Z$.  
    By \cite[Theorem~2.39]{birkar2026steindegreefieldscharacteristic}, we can write
    \[ K_{X'} + B' + \mb{M}_{X'} \sim_{\RR} g^*(K_{Y'}+C' + \mb{N}_{Y'}), \]
    where $(Y', C' + \mb{N})$ is a g-lc g-pair.  
    By \cite[Lemmas~2.49, 2.51]{birkar2026steindegreefieldscharacteristic},
    $X'\to Z$ and $X'\to Y'$ are contractions of Fano type.
    Hence by \cite[Theorem~1.2]{birkar2026steindegreefieldscharacteristic},
    $\sdeg ((S')^{\nor}/T')$ is bounded from above depending only on $d,t$.
    
    By \cite[Lemma~2.38]{birkar2026steindegreefieldscharacteristic}, we get
    \[ \coeff_{T'} C' \ge t. \]
    Moreover, by construction, it is evident that
    \[ K_{Y'}+C' + \mb{N}_{Y'} \sim_{\RR}0/Z. \]
    Thus, $(Y', C' + \mb{N})\to Z$ is a generalised log
    Calabi-Yau fibration over $k$.  Additionally, by Lemma~\ref{birational-fano-type},
    $Y'\to Z$ is a birational contraction of Fano type.
    Since we assume Theorem~\ref{dim-d-t} in dimension $d$, 
    we deduce that $\sdeg ((T')^{\nor}/T)$ is bounded from above depending only on $d,t$.
    Therefore, by \cite[Lemma~2.9]{birkar2026steindegreefieldscharacteristic}, we can conclude that
    \[ \sdeg (S^{\nor}/T) = \sdeg ((S')^{\nor}/T) \le \sdeg ((S')^{\nor}/T')\cdot \sdeg ((T')^{\nor}/T), \]
    which is bounded from above depending only on $d,t$.
\end{proof}


\begin{proof}[Proof of Theorems~\ref{dim-d-t} and~\ref{thm-vertical-introduction}]
    We proceed by induction on $d$.    
    Let $(X, B + \mb{M})\to Z$ be a Fano type generalised log Calabi-Yau fibration
    over a field $k$ of characteristic zero,
    and let $S$ be an irreducible component of $B$ with $\coeff_S B \ge t$.
    Assume $\dim X = 1$.  In this case, $S$ is a closed point on $X$.
    If $\dim Z = 0$, then Theorem~\ref{thm-vertical-introduction}
    follows from \cite[Theorem~1.1]{birkar2026steindegreefieldscharacteristic}.
    If $\dim Z = 1$, then $X\to Z$ is an isomorphism of nonsingular curves,
    so $S$ is isomorphic to its image on $Z$.  Hence the Stein degree is equal to 1.
    Thus, Theorem~\ref{thm-vertical-introduction} holds for $d=1$ and each $t\in (0,1]$.
    Then by Lemmas~\ref{vertical-to-birational} and~\ref{birational-to-vertical},
    we can conclude that Theorems~\ref{dim-d-t} and~\ref{thm-vertical-introduction}
    hold for every $d\in \NN$ and every $t\in (0, 1]$.
\end{proof}


\section{\bf Examples}\label{s-examples}

In this section, we construct examples showing that the dimension
dependence in our bounds cannot be omitted or replaced by dependence
only on the relative dimension and the coefficient lower bound.


\subsection{Dependence on dimension}\label{exa-dim}
In this subsection, we give examples showing that 
the dependence on dimension in Theorem~\ref{dim-d-t} cannot be removed.

\medskip
(1)
Let $n$ be a natural number $\ge 2$.
Let $T$ be a smooth integral scheme over a field $k$ of characteristic zero, and let 
\[ p\colon U\to T \] 
be a finite \'etale morphism of degree $n$ with $U$ integral.
(For example, let $k \coloneqq \CC$, and let
\[ T\coloneqq \spec\CC [t,t^{-1}]~~\mbox{and}~~ U\coloneqq \spec\CC [u,u^{-1}]. \]
Then define the morphism $p\colon U\to T$ by $t\mapsto u^n$.)
Set $\mcA \coloneqq p_*\mcO_U$, which is a locally free $\mcO_T$-algebra of rank $n$.  Let
\[ Z\coloneqq \ul{\spec}_T \big(\sym_{\mcO_T} (\mcA^{\vee})\big). \]
In other words, $Z$ is the total space of the vector bundle $\mcA$.
Let $T_0\subset Z$ be the zero section of $\pi\colon Z\to T$.  Then $T_0\cong T$ via $Z\to T$.

Let $\xi$ be the tautological section of $\pi^*\mcA$.
Let $D\subset Z$ be the reduced divisor defined by the vanishing of the norm of $\xi$;
see \cite[(6.5.1)]{EGA-II}.  If $V\to T$ is an \'etale morphism
such that $V\times_T U$ splits into disjoint copies of $V$, then the 
pullback of $D$ to $V$ is an snc divisor.  
In particular, $(Z, D)$ is an lc pair.
Let $Z_1\to Z$
be the blow-up along the zero section $T_0\subset Z$,
and let $E_1$ be the exceptional divisor.  Then we get
\[ E_1 = \PP_{T_0}(\mcA^{\vee}). \]
In other words, $E_1\to T_0$ is the projective bundle of lines in $\mcA$.

\medskip
(2)
We claim that there exists a closed immersion
\[ \iota\colon U\hookrightarrow E_1 \]
over $T_0$, where we view $U$ as a $T_0$-scheme via the isomorphism $T\cong T_0$.
As $p$ is finite \'etale, the diagonal $\Delta_U$ of $U\to U\times_T U$ is open and closed in $U\times_T U$.  Write
\[ U\times_T U = \Delta_U \amalg W. \]
Let $e$ be the regular function on $U\times_T U$ that takes value 1 on $\Delta_U$ and 0
on $W$.  Let 
\[ \projection_1 \colon U\times_T U \to U \]
be the canonical projection to its first factor.  
By \cite[Corollaire~(1.5.2)]{EGA-II}, we get
\[ p^*\mcA \cong (\projection_1)_*(\mcO_{U\times_T U}). \]
Then $e$ generates a direct summand $\mcO_U e\subset p^* \mcA$ of rank one.
Indeed, $\projection_1|_{\Delta_U}$ is an isomorphism, and $\projection_1|_{W}$ is 
a finite \'etale morphism of degree $n-1$.  Hence we can write
\[ p^*\mcA = \mcO_U\oplus \mcB, \]
where $\mcB \coloneqq (\projection_1|_W)_* \mcO_W$ is a locally free sheaf of rank $n-1$ on $U$.
Then $e$ is $(1,0)$ in the decomposition of $p^*\mcA$.
This shows that $e$ generates a direct summand of $p^*\mcA$.

Recall that $E_1=\PP_{T_0}(\mcA^\vee)$.  Then
\[ E_1\times_{T_0} U = \PP_U(p^*\mcA^\vee). \]
Dualising the split inclusion $\mcO_U e \hookrightarrow p^* \mcA$ gives a surjection
\[ p^*\mcA^{\vee} \to (\mcO_U e)^{\vee}\cong \mcO_U. \]
By \cite[Th\'eor\`eme~(4.2.4)]{EGA-II}, this corresponds to a section
\[ \sigma\colon U\to \PP_U(p^*\mcA^{\vee}) \]
of the projection $\PP_U(p^*\mcA^{\vee}) \to U$.  Denote by $\iota\colon U\to E_1$
the composite morphism
\[ U\xrightarrow{\sigma} E_1\times_{T_0} U \xrightarrow{\projection_{E_1}} E_1. \]
It is evident that $\iota$ is a morphism over $T_0$.
    
Now we show that $\iota$ is a closed immersion.
Pick a finite \'etale cover $T'\to T$ that splits $p$, i.e.
\[ U'\coloneqq U\times_T T' = \coprod_{i=1}^n T_i', \]
where $T_i'\cong T'$ for each $1\le i \le n$.
By \cite[Proposition~III.9.3]{Hart}, we get $\mcA|_{T'} \cong \mcO_{T'}^{\oplus n}$.
Denote by $e_1, \dots, e_n$ the standard idempotents of $\mcA|_{T'}$.  Then
\[ E_1'\coloneqq E_1\times_{T_0} T' \cong \PP_{T'}^{n-1}, \]
where again we identify $T_0$ and $T$ with the isomorphism $T_0\cong T$ via $\pi\colon Z\to T$.  Let
\[ \iota'\colon U'\to E_1' \]
be the base change of $\iota$ via $T'\to T$.  The scheme $U'\times_{T'} U'$
is the disjoint union of copies $T_i'\times_{T'} T_j' \cong T'$ indexed by pairs $(i, j)$
with $1\le i, j \le n$.  The diagonal $\Delta_{U'} \subset U'\times_{T'} U'$
consists exactly of the copies with $i=j$.
Denote by
\[ \projection_1' \colon U'\times_{T'} U' \to U' \]
the canonical projection to the first factor.  
Then as a section of $(\projection_1')_*(\mcO_{U'\times_{T'}U'})$, 
the pullback $e'$ of $e$ to $U'$ takes value 1 on the $(i,i)$-copy
and 0 on every $(i,j)$-copy with $i\neq j$.
Hence on the $i$-th component $T_i'$ of $U'$, $e'$ is equal to $e_i$.

Denote by $p'\colon U'\to T'$ the base change of $p$.  Let
\[ \sigma'\colon U'\to E_1'\times_{T'} U' \]
be the base change of $\sigma$. By construction,
\[ \iota' = \projection_{E_1'}\circ\sigma'. \]
Moreover, $\sigma'$ corresponds to the quotient obtained by dualising
$\mcO_{U'}e'\hookrightarrow p'^*(\mcA|_{T'})$.
Since $e'|_{T_i'}=e_i$, the restriction of this quotient to $T_i'$ is
\[ \bigoplus_{j=1}^n\mcO_{T_i'}e_j^\vee
    \to \mcO_{T_i'},~\sum_{j=1}^n a_j e_j^\vee\mapsto a_i. \]
Thus, under the identification
\[ (E_1'\times_{T'}U')|_{T_i'} \cong \PP_{T_i'}^{n-1}, \]
the restriction $\sigma'|_{T_i'}$ is the $i$-th coordinate section.
Note that the projection
\[ E_1'\times_{T'}T_i'\to E_1' \]
is an isomorphism as $T_i'\to T'$ is an isomorphism.
Therefore, composing $\sigma'|_{T_i'}$ with this projection
shows that $\iota'|_{T_i'}$ is also the $i$-th coordinate section of $E_1'\cong\PP_{T'}^{n-1}$.
Hence $\iota'\colon U' \to E_1'$ is a closed immersion over $T'$.
By \cite[Proposition~(2.7.1)]{EGA-IV-2},
we can conclude that $\iota\colon U\to E_1$ is a closed immersion over $T$.

\medskip
(3)
Let $X\to Z_1$ be the blow-up along $\iota(U)$, and let $S$ be the exceptional divisor.
Denote by $f\colon X\to Z$ the composite $X\to Z_1 \to Z$.
Let $B$ be the birational transform of $D$ to $X$.  Write
\[ B_1\coloneqq D_1 + E_1 ~~\mbox{and}~~ B \coloneqq B_1^{\sim} + S, \]
where $D_1$ is the birational transform of $D$ to $Z_1$,
and $B_1^{\sim}$ is the birational transform of $B_1$ to $X$.
After passing to a finite \'etale cover of $T$ splitting $p\colon U\to T$, we get
\[ K_X + B = f^*(K_Z + D). \]
In particular, $(X, B)$ is lc and $\coeff_S B = 1$.
Moreover, by \cite[Lemma~2.50]{birkar2026steindegreefieldscharacteristic}, 
$X\to Z$ is of Fano type.
By construction, $S\to \iota(U)$ is a $\PP^{n-1}$-bundle.
It is evident that the image of $S$ in $Z$ is $T_0$.  Thus, we get
\[ \sdeg (S^{\nor}/T_0) = \sdeg (S/T_0) = n. \]
As $\dim X = \dim T + n$, we can conclude that
the dependence on dimension in Theorem~\ref{dim-d-t} can not be omitted.


\subsection{Insufficiency of relative dimension}\label{exa-rel-dim}
Let $n$ be a natural number $\ge 2$.
Let $T$ be a smooth integral scheme over a field $k$ of characteristic zero, and let 
\[ p\colon U\to T \] 
be a finite \'etale morphism of degree $n$ with $U$ integral.  Let
\[ f \colon (X, B) \to Z \]
be the birational log Calabi-Yau fibration constructed in \S\ref{exa-dim},
and let $S$ be the nonsingular divisor on $X$ with $\sdeg (S/T_0)= n$,
where $T_0$ denotes the image of $S$ in $Z$.  Let
\[ Y \coloneqq X\times_{\spec k} \PP_k^1. \]
Denote by 
\[ \projection_1\colon Y\to X~~\mbox{and}~~\projection_2 \colon Y\to \PP_k^1 \]
the projections to the first and second factors.  Let
\[ \Delta\coloneqq \projection_1^* B + \projection_2^*(\set{0} + \set{\infty}). \]
Let
\[ g\coloneqq f\circ\projection_1\colon Y\to Z. \]
It is evident that $g$ is a contraction, $(Y, \Delta)$ is lc, and 
\[ K_Y + \Delta = \projection_1^*(K_X+B) + \projection_2^*(K_{\PP_k^1} + \set{0} + \set{\infty}) \sim_{\QQ} g^*(K_Z + D), \]
where $D$ is the reduced divisor on $Z$ defined in \S\ref{exa-dim}(1).
Hence $(Y, \Delta)\to Z$ is a log Calabi-Yau fibration of relative dimension one.

As $f\colon X\to Z$ is of Fano type, pick a $\QQ$-divisor $\Gamma$ on $X$ such that 
$(X, \Gamma)$ is klt, and $A\coloneqq -(K_X + \Gamma)$ is ample$/Z$.  Then $(Y, \projection_1^*\Gamma)$
is also klt, and
\[ -(K_Y + \projection_1^*\Gamma) = \projection_1^*A + \projection_2^*(-K_{\PP_k^1}) \]
is ample$/Z$.  Thus, $Y\to Z$ is of Fano type.

Let
\[ \wt{S} \coloneqq S\times_{\spec k}\PP_k^1. \]
Then $\wt{S}$ has coefficient 1 in $\Delta$, and
\[ \sdeg (\wt{S}/T_0) = n. \]
The Stein degrees of vertical divisors are therefore unbounded even when the relative   
dimension is one and the boundary coefficient is fixed at one.


\medskip

\printbibliography[title={\bf References}]

\vspace{1em}
 
\noindent\small{Xintong Jiang} 

\noindent\small{\textsc{Qiuzhen College, Tsinghua University, Beijing, China} }

\noindent\small{Email: \texttt{xt-jiang21@mails.tsinghua.edu.cn}}

\vspace{1em}
 
\noindent\small{Santai Qu} 

\noindent\small{\textsc{Beijing Institute of Mathematical Sciences and Applications, Beijing, China} }

\noindent\small{Email: \texttt{qusantai@bimsa.cn}}

\end{document}